\documentclass[11pt,a4paper]{amsart}
\usepackage{
    amsmath,  amsfonts, amssymb,  amsthm,   amscd,
    gensymb,  graphicx, comment,  etoolbox, url,
    booktabs, stackrel, mathtools,enumitem, mathdots,  microtype, lmodern,    mathrsfs, graphicx, tikz,  longtable,tabularx, float, tikz, pst-node, tikz-cd, multirow, tabularx, amscd,  bm, array, makecell, diagbox, booktabs,ragged2e
}
\usepackage{xcolor}
\usepackage[utf8]{inputenc}
\usepackage{microtype, fullpage, wrapfig,textcomp,mathrsfs,csquotes,fbb}
\usepackage[colorlinks=true, linkcolor=blue, citecolor=blue, urlcolor=blue, breaklinks=true]{hyperref}
\usepackage[capitalise]{cleveref}
\usepackage{todonotes}
\usepackage{tikz-cd}
\usepackage[foot]{amsaddr}
\usetikzlibrary{arrows}

\newtheorem{theorem}{Theorem}[section]

\newtheorem{cor}[theorem] {Corollary}
\newtheorem{definition}[theorem]{Definition}
\newtheorem{example}[theorem]{Example}
\newtheorem{lemma}[theorem]{Lemma}
\newtheorem{notation}[theorem]{Notation}

\newtheorem{remark}[theorem]{Remark}
\begin{document}
\title{Band of topological groups}

\author{S\lowercase{unil} K\lowercase{umar} M\lowercase{aity}$^1$ \lowercase{and} M\lowercase{onika} P\lowercase{aul}$^{2, \ast}$}

\address{$^1$University of Calcutta}
\email{skmpm@caluniv.ac.in}
\address{$^2$Indian Institute of Technology Madras}
\address{$^{\ast}$Corresponding author}
\email{paulmonika007@gmail.com}

\begin{abstract} 
In this article, we construct a band of topological groups from a cryptogroup. Also, we prove that a band of topological groups is metrizable if and only if each $\mathcal{H}$-class is metrizable. Finally, we demonstrate that if $S$ is a band of topological groups and $N$ is a full normal subcryptogroup of $S$, then $S/N$ is Hausdorff if and only if $\rho_{_N}$ is closed in $S \times S$ if and,  $\rho_{_N}$ is closed in $S \times S$ if and only if $N$ is closed in $S$. 
\end{abstract}
\maketitle

\section{Introduction}
In a semigroup $S$, an element $a\in S$ is said to be an idempotent element if $a^2 = a$, and the set of all idempotent elements in $S$ is denoted by $E(S)$. A semigroup $S$ is said to be a band if $E(S) = S$. An element $a$ in a semigroup $S$ is said to have an inverse $x$ in $S$ if $xax = x$ and $axa = a$, and the set of all inverse elements of $a \in S$, if it exists, is denoted by $V(a)$. Moreover, in a semigroup $S$, an element $a$ is said to be regular (algebraic property) if $axa = a$, for some $x \in S$ and in this case, if we take $y = xax$, then $aya = a$, $yay = y$ and so $V(a) \neq \emptyset$. A regular (algebraic property) semigroup is a semigroup in which each element is regular (algebraic property). In a regular (algebraic property) semigroup $S$, the Green's relations \cite{how} $\mathcal{L}, \mathcal{R}$ and $\mathcal{H}$ are defined by:  for any $a,b \in S$,

\begin{center}
$a \, \mathcal{L} \, b$ if and only if $Sa = Sb$;\\
$a \, \mathcal{R} \,  b$ if and only if $aS = bS$;\\
$\mathcal{H} = \mathcal{L} \cap \mathcal{R}$.
\end{center}

An element $a$ in a semigroup $S$ is called completely regular \cite{pr} (algebraic property) if $a = axa$ and $ax = xa$, for some $x \in S$. A semigroup $S$ is said to be a completely regular semigroup \cite{pr} (algebraic property) if each of its elements is completely regular (algebraic property). A semigroup $S$ is said to be cryptic if $\mathcal{H}$ is a congruence on $S$. A cryptic completely regular semigroup is called a cryptogroup. In a cryptogroup $S$, for each $x \in S$, we denote the $\mathcal{H}$-class containing $x$ by $H_{_x}$ and the identity element of the group $H_{_x}$ is denoted by $x^0$. In addition, in a cryptogroup $S$, for each $x \in S$, there exists a unique $y \in H_{_x}$ such that $xy = yx = x^0$ and we denote this $y$ by $x^{-1}$. A nonempty subset $K$ of a cryptogroup $S$ is said to be a full normal subcryptogroup of $S$ if

$(i)$ $K$ is full, i.e., $E(S) \subseteq K$;

$(ii)$ $K$ is a subcryptogroup of $S$, i.e., $xy, x^{-1} \in K$, for all $x, y \in K$;

$(iii)$ $K$ is normal, i.e., $sks^{-1} \in K$, for all $s \in S$ and $k \in K$.

 A congruence $\rho$ on a semigroup $S$ is said to be a band congruence if $S/ \rho$ is a band. A semigroup $S$ is said to be a band $B$ of groups $G_{_\alpha}(\alpha \in B)$ if $S$ admits a band congruence $\rho$ on $S$ such that $S/\rho = B$ and each $G_{_\alpha}$ is a $\rho$-class mapped onto $\alpha$ by the natural epimorphism ${\rho}^\# : S \longrightarrow B$. It is interesting to note that $\mathcal{H}$ is a band congruence on a cryptogroup $S$ and each $\mathcal{H}$-class is a group. Hence a semigroup $S$ is a cryptogroup if and only if $S$ is a band of groups.

Let $X$ and $Y$ be topological spaces. Any open set $A$ in $X$ which is also closed is a clopen set.  A mapping $f: X \longrightarrow Y$ is said to be continuous \cite{mun} if for any open set $V$ in $Y$, $f^{-1}(V)$ is open in $X$. Also, a mapping $f: X \longrightarrow Y$ is said to be open if for any open set $U$ in $X$, $f(U)$ is open in $Y$. A continuous bijective mapping $f: X \longrightarrow Y$ which is open is a homeomorphism. A continuous surjective mapping $g: X \longrightarrow Y$ is said to be a quotient mapping if for any subset $V$ of $Y$, $g^{-1}(V)$ is open in $X$ implies that $V$ is open in $Y$. 

A topological space $X$ is said to be a

\begin{itemize}
\item[-] $T_0$ space if for any two distinct points $x,y \in X$, there is an open set $U$ which contains exactly one of $x$ and $y$;
\item[-] $T_1$ space if for any two distinct points $x,y \in X$, there exist open sets $U,V$ such that $x\in U, y \notin U$ and $x \notin V, y \in V$;
\item[-] $T_2$ or Hausdorff space if for any two distinct points $x,y \in X$, there exist two disjoint open sets $U$ and $V$ containing $x$ and $y$ respectively in $X$;
\item[-] regular (topological property) if for any $x\in X$ and any closed set $F$ with $x \notin F$, there exist two disjoint open sets $U$ and $V$ containing $x$ and $F$ respectively in $X$;
\item[-] completely regular (topological property) if for any $x\in X$ and any closed set $F$ with $x \notin F$, there is a continuous mapping $f : X \longrightarrow [0,1]$ such that $f(x) = 1$ and $f(F) = \{0\}$.
\item[-] normal (topological property) if for any two disjoint closed sets $A$ and $B$ in $X$, there exist two disjoint open sets $U$ and $V$ containing $A$ and $B$ respectively in $X$.
\end{itemize}

A family of subsets $\mathcal{F}$ of a topological space $X$ is said to be discrete if, for each $x \in X$, there exists an open set $U$ containing $x$ such that $U$ intersects at most one element in $\mathcal{F}$. A topological space $X$ which cannot be expressed as a union of two nonempty disjoint open sets is called a connected space. A topological space $X$ is said to be locally connected at $x\in X$ if any open set containing $x$ contains a connected open set containing $x$. A topological space $X$ which is locally connected at each of its points is called a locally connected space. A topological space $X$ is said to be compact if every open cover of $X$ has a finite subcover. Moreover, a topological space is locally compact at $x \in X$ if there is an open set $W$ containing $x$ such that $\overline{W}$ is compact. A topological space that is locally compact at each $x \in X$ is called a locally compact space. It is worth mentioning that a continuous image of a compact space is also a compact space. Any closed subset $A$ of a compact space $X$ is also compact and any compact subset of a Hausdorff space is closed. In addition, any closed subset of a locally compact Hausdorff space $X$ is locally compact.
 
A semigroup $S$ endowed with a topology $\tau$ is said to be a topological semigroup if the mapping $\mu: S\times S \longrightarrow S$ defined by $(x,y) \mapsto xy$, is continuous, i.e., for any $x,y \in S$ and any open set $W$ containing $xy$ in $S$, there exist open sets $U$ and $V$ containing $x$ and $y$ respectively in $S$ such that $UV \subseteq W$. A topological semigroup $(S, \tau)$ is said to be a  

\begin{itemize}
\item[-] topological cryptogroup if $S$ is a cryptogroup and the mapping $\gamma: S \longrightarrow S$ defined by $x  \, \mapsto \, x^{-1}$, is continuous. 

\item[-] topological band if $S$ is a band.

\item[-] band $B$ of topological groups $(G_{_\alpha}, \tau_{_\alpha})(\alpha \in B)$ if $S$ admits a band congruence 
$\rho$ such that $S/ \rho = B$, each $G_{_\alpha}$ is a $\rho$-class mapped onto $\alpha$ by the natural semigroup epimorphism $\rho^{\#} : S \longrightarrow B$ and $\bigcup\limits_{\alpha \in B}\tau_{_\alpha}$ forms a base for the topology $\tau$.
\end{itemize}
A relation $\rho$ on a band of topological groups $S$ is said to be closed if $\rho$ is a closed subset of $S \times S$.

In Section 2, we establish the necessary and sufficient conditions for a topological cryptogroup to be a band of topological groups. The construction of a band of topological groups from a cryptogroup has been studied in Section 3. In Section 4, some topological properties have been established and Section 5 is devoted to studying the full subcryptogroup of a cryptogroup. Finally, in Section 6, we study the full normal subcryptogroup of a topological cryptogroup. 

\section{Band of topological groups}
In this section, we first give some examples of bands of topological groups.
\begin{example}
    We consider the semigroup $S = (\mathbb{Z}_{10}, \cdot_{10})$. Then $S$ is a topological semigroup with respect to the topology $\tau$ generated by the sets $\{\bar{0}\}, \{\bar{5}\}, \{\bar{1}\}, \{\bar{9}\}, \{\bar{3}\}, \{\bar{7}\}, \{\bar{2}, \bar{4}, \bar{6}, \bar{8}\}$. It can be easily verified that $(S, \tau )$ is a band of topological group.
\end{example}
\begin{example}
    Let us consider the semigroup $S = (\mathbb{Z}_{15}, \cdot_{15})$. Then $S$ is a topological semigroup with respect to the topology $\tau$ generated by the sets $\{\bar{0}\},  \{\bar{3}\},  \{\bar{5}\},\{\bar{6}\}, \{\bar{9}\}, \{\bar{10}\},  \{\bar{12}\}, \{\bar{1}, \bar{2}, \bar{4}, \bar{7}, \bar{8}, \bar{11}, \bar{13}, \bar{14}\}$. Then $(S, \tau )$ is a band of topological group.
\end{example}

From the above two examples, it is clear semigroups mentioned in the above examples are topological cryptogroups, which assert a connection between topological cryptogroups and bands of topological groups. A natural question may arise: is any topological cryptogroup a band of topological groups? By a counterexample, we show the false assumption, that is, a topological cryptogroup may not be a band of topological groups.

\begin{example}
	Consider the topological semigroup  $(S, \cdot_6, \tau)$, where $S$ is the set of all integers modulo $6$ and $\tau = \{S, \emptyset, \{\bar{0}\},  \{\bar{3}\}, \{\bar{1}, \bar{2}, \bar{4}, \bar{5}\}\}$. Then $(S, \tau)$ is a topological cryptogroup which is not a band of topological groups.
\end{example}

We now establish a necessary and sufficient condition on a topological cryptogroup $S$ to be a band of topological groups. For this purpose, we register a topological property.

\begin{lemma}\label{top. prop.}
For any equivalence relation $\rho$ on a topological space $S$, each $\rho$-class is open in $S$ if and only if for each $x \in S$ and each open set $G$ containing $x$, there exists an open set $U$ in $S$ such that $x \in U \subseteq G \cap x \rho$, where $x\rho$ is the $\rho$-class containing the element $x$.
\end{lemma}

\begin{proof}Let each $\rho$ class be open in $S$. Let $x \in S$ and $G$ be an open set containing $x$. Then $G \cap x\rho$ is open in $S$ and set $U = G \cap x \rho$. Then we have $x \in U \subseteq G \cap x \rho$. 

Conversely, let for each $x \in S$ and each open set $G$ containing $x$, there exists an open set $U$ in $S$ such that $x \in U \subseteq G \cap x \rho$. Let $H$ be a $\rho$-class and $x \in H$. Now $S$ being open with $x \in S$ there exists an open set $U$ such that $x \in U \subseteq S \cap H = H$ and hence $H$ is open in $S$.
\end{proof}

\begin{theorem}\label{main theorem}
A topological semigroup $(S, \tau)$ is a band of topological groups if and only if $S$ is a topological cryptogroup satisfying the property that for each $x \in S$ and each open set $G$ in $S$ containing $x$ there is an open set $U$ in $S$ such that $x\in U \subseteq G \cap H_{_x}$.
\end{theorem}

\begin{proof}
Let $(S, \tau)$ be a band $B$ of topological groups $(G_{_\alpha}, \tau_{_\alpha} )$. Then $S$ is a cryptogroup. To show the continuity of the mapping $\gamma$, let $x\in S$ and $W$ be an open set in $S$ containing $x^{-1}$. Then there is a basic open set $U$ in $S$ such that $x^{-1} \in U \subseteq W$. This implies $U \in \tau_{_\alpha}$, for some $\alpha \in B$ with $x^{-1} \in G_{_\alpha}$. Now $G_{_\alpha}$ being a topological group, $U^{-1}$ is open in $G_{_\alpha}$ and hence in $S$. Moreover, $x \in U^{-1}$ and $(U^{-1})^{-1} \subseteq W$. It follows that $(S, \tau)$ is a topological cryptogroup. Now, let $y \in S$ and $G$ be an open set containing $y$ in $S$. Since, $\cup \tau_{_\alpha}$ generates $\tau$, there is a basic open set $V$ in $\cup \tau_{_\alpha}$ such that $y \in V \subseteq G$. Now, $V \in \tau_{_\beta}$, for some $\beta \in B$ with $y \in G_{_\beta} = H_{_y}$. Hence $y \in V \subseteq G \cap H_{_y}$. 

Conversely, let $(S, \tau)$ be a topological cryptogroup satisfying the property that for each $x \in S$ and each open set $G$ in $S$ containing $x$ there is an open set $U$ in $S$ such that $x\in U \subseteq G \cap H_{_x}$. Then $\mathcal{H}$ is a band congruence and each $\mathcal{H}$-class is a topological group. Let $B = S/\mathcal{H}$. Then each $\mathcal{H}$-class mapped onto $\alpha \, (\in B)$ by the natural epimorphism $\rho^{\#} : S \longrightarrow B$.  By Lemma \ref{top. prop.}, each $\mathcal{H}$-class is open in $S$. Let  $\tau_{_x}$ be the topology on $H_{_x}$, for each $x \in X$. To complete the proof, it remains to show that $\cup \tau_{_x}$ generates $(S, \tau)$. Let $\tau'$ be a topology generated by $\cup \tau_{_x}$. Clearly, each element in $\tau_{_x}$ is open in $H_{_x}$, for each $x\in S$ and so open in $(S, \tau)$. Therefore, $\cup \tau_{_x} \subseteq \tau$ and thus $\tau' \subseteq \tau$. For the reverse inclusion, let $W \in \tau$. Then for each $z \in W$, $W\cap H_{_z}$ is open in $H_{_z}$ and so $W \cap H_{_z} \in \tau^\prime$. Therefore, $W = \displaystyle{\bigcup_{z \in W}(W\cap H_z} ) \in \tau'$ and so $\tau \subseteq \tau^\prime$. Consequently, $\cup\tau_\alpha$ generates $\tau$
\end{proof}

\begin{cor}
Each $\mathcal{H}$-class in a band of topological groups is a clopen topological group. 
\end{cor}

\begin{proof}
Let $S$ be a band of topological groups. Clearly, each $\mathcal{H}$-class is a topological group. Moreover, by Lemma \ref{top. prop.} and Theorem \ref{main theorem}, it follows that each $\mathcal{H}$-class is open in $S$. Again, since $S$ is a cryptogroup, so $S$ is union of its $\mathcal{H}$-classes and hence for each $\mathcal{H}$-class $H$, $H = S \setminus \displaystyle{\bigcup_{x \in S\setminus H} H_x}$. Since  $\displaystyle{\bigcup_{x \in S\setminus H} H_x}$ is open in $S$, it follows that its complement $H = S\setminus \displaystyle{\bigcup_{x \in S\setminus H} H_x}$ is closed in $S$, i.e., $H$ is clopen in $S$. 
\end{proof}

Combining Lemma \ref{top. prop.} and Theorem \ref{main theorem}, we have the following corollary.

\begin{cor}\label{main theorem 2}
A topological cryptogroup  $S$ is a band of topological groups if and only if each $\mathcal{H}$-class is open in $S$.
\end{cor}

Let $S$ be a cryptogroup endowed with a topology $\tau$. Consider the natural epimorphism $\pi: S \longrightarrow S/ \mathcal{H}$. Now we define a topology on $S/\mathcal{H}$ as follows : a subset $U$ of $S/ \mathcal{H}$ is open in $S/ \mathcal{H}$ if and only if $\pi^{-1}(U)$ is open in $S$. It is immediate that in this case, $\pi$ is a quotient map.

\begin{theorem}
Let $S$ be a cryptogroup. Then $S/ \mathcal{H}$ is a discrete topological band if and only if $S$ admits a topology $\tau$ for which $(S,\tau)$ is a band of topological groups.
\end{theorem}

\begin{proof}
Let $\mu$ and $\mu_1$ be the corresponding binary operations defined on $S$ and $S/ \mathcal {H}$ respectively. Let $S/ \mathcal{H}$ be a discrete topological band. Consider the natural epimorphism $\pi : S \longrightarrow S/\mathcal{H}$ defined by $\pi (x) = H_{_x}$. Define a topology  $\tau$ on $S$ as follows : $\tau = \{\pi^{-1}(U): U\,\, \text{is open in}\,\, S/ \mathcal{H} \}$. With the topology $\tau$ on $S$, $\pi$ is continuous. Moreover, $\mu_1 \circ (\pi \times \pi) = \pi \circ \mu$. To show the continuity of $\mu$, let $\pi^{-1}(W)$ be an open set in $S$, where $W$ is open in $S/\mathcal{H}$. Then $(\pi \times \pi)^{-1} (\mu_1^{-1}(W)) = \mu^{-1}(\pi^{-1}(W)) $. Now, $S/\mathcal{H}$ being a topological band, it follows that $\mu^{-1}(\pi^{-1}(W))$ is open in $S\times S$. Hence $\mu$ is continuous. For the continuity of $\gamma$, let $\pi^{-1}(W_1)$ be an open set in $S$. It can be easily shown that $\gamma^{-1}(\pi^{-1}(W_1)) = \pi^{-1}(W_1)$. So, $\gamma$ is continuous and hence $S$ is a topological cryptogroup. By the construction of $\tau$, it is immediate that each $\mathcal{H}$-class is open in $S$. So by Corollary \ref{main theorem 2}, it follows that $S$ is a band of topological groups. 

Conversely, let $S$ be a band of topological groups. Since each $\mathcal{H}$-class is open in $S$, $S/\mathcal{H}$ is a discrete space which is a band. The topology on $S/ \mathcal{H}$ being discrete, the binary operation $\mu_1$ is continuous, and hence $S/ \mathcal{H}$ is a discrete topological band.
\end{proof}

Now, we establish some sufficient conditions for a cryptogroup for being a band of topological groups.

\begin{theorem}
Let $S$ be a topological semigroup such that $S$ is a cryptogroup. Also, assume that for each $e \in E(S)$ and for any open neighborhood $U$ of $e$, there exists an open neighborhood $V$ of $e$ such that $V^{-1} \subseteq U$ and each $\mathcal{H}$- class is open in $S$. Then $S$ is a band of topological groups.
\end{theorem}

\begin{proof}
At first we prove that if $G$ be an open set in $S$ and $x \in S$, then $(xG)^*$ and $(Gx)^*$ are open in $S$, where $(xG)^* = \{y \in S : x^{-1}y \in G \, \mbox{and} \, x^0 = y^0\}$ and $(Gx)^* = \{y \in S : yx^{-1} \in G \, \mbox{and} \, x^0 = y^0\}$. It follows from the fact that $S$ is a topological semigroup and each $\mathcal{H}$ - class is open in $S$. Now we show that the mapping $x \mapsto x^{-1}$ is continuous. Let $W$ be an open set in $S$ and $x \in S$ with $x^{-1} \in W$. This implies that $x^0x^{-1} \in W$. Then there exists an open set $U$ containing $x^0$ in $S$ such that $Ux^{-1} \subseteq W$. Now by the given assumption, there exists an open set $V$ containing $x^0$ in $S$ such that $V^{-1} \subseteq U$. This implies that $V^{-1}x^{-1} \subseteq Ux^{-1} \subseteq W$. Then $(V^{-1}x^{-1})^* \subseteq (Ux^{-1})^* \subseteq (Ux^{-1}) \subseteq W$. Now $(V^{-1}x^{-1})^* = ((xV)^*)^{-1}$. Thus, we have $(xV)^*$ is an open set containing $x$ with $((xV)^*)^{-1} \subseteq W$. Therefore, the mapping $x \mapsto x^{-1}$ is continuous and hence $S$  is a band of topological groups.
\end{proof}

\section{Construction of band of topological groups}
In this section, we construct a band of topological groups from a cryptogroup. At first, we establish some results on open sets and neighbourhood systems at an idempotent element of a band of topological groups. For this purpose, let us define the following definition.

\begin{definition}
Let $S$ be a cryptogroup and $U, V$ be two nonempty subsets of $S$. Define a set $(UV)^* = \bigcup \{(uV)^* : u \in U\}$, where $(uV)^* = \{x\in S : u^{-1}x \in V \,\, \text{and}\,\, u^0 = x^0\}$. It can be easily verify that $(uV)^* \subseteq uV$ and $(UV)^* = \bigcup \{(Uv)^* : v \in V\}$, where $(Uv)^* = \{x \in S : xv^{-1} \in U \,\, \text{and}\,\, x^0 = v^0\}$. Moreover, $((xU)^*y)^* = (x(Uy)^*)^* = (xUy)^*$, where $(xUy)^* = \{u\in U : x^{-1}uy^{-1} \in U\, \, \text{and} \,\, x^0 = u^0 = y^0\}$. In this connection, it is interesting to note that if two elements $x, y \in S$ are not $\mathcal{H}$-related, then $(xUy)^* = \emptyset$.
\end{definition}

\begin{theorem}
Let $D$ be a dense subset of a band of topological groups $S$. Then for any open set $U$ containing $E(S)$, $(UD)^* = S= (DU)^*$.
\end{theorem}

\begin{proof}
Clearly, $(UD)^* \subseteq S$. Let $s\in S$. Then $(s^{-1}U)^* = \{x\in S : sx \in U \, \, \text{and} \,\, s^0 = x^0\}$ is an open set containing $s^{-1}$. It follows that $((s^{-1}U)^*)^{-1}$ is an open set containing $s$. Since $\overline{D} = S$, $((s^{-1}U)^*)^{-1} \cap D \neq \emptyset$. Let $t \in ((s^{-1}U)^*)^{-1} \cap D$. Then $t^{-1} \in (s^{-1}U)^*$. This implies that $st^{-1} \in U$ and  $s^0 = (t^{-1})^0$. Set $u =st^{-1}$. Then $u^{-1}s = ts^{-1}s = t \in D$ and $u^0 = (st^{-1})^0 = (s^0(t^{-1})^0)^0 = s^0$. It follows that $s \in (uD)^* \subseteq (UD)^*$. Therefore, $(UD)^* = S$. Similarly, one can easily show that $(DU)^* = S$. 
\end{proof}

\begin{theorem}
Let $S$ be a band of topological groups. Then for every subset $A$ of $S$ and every open set $U$ containing $E(S)$, $\overline{A} \subseteq (AU)^*$.
\end{theorem}

\begin{proof}
Let $A$ be a subset of $S$ and $U$ be an open set containing $E(S)$. For each $e \in E(S)$, there exists an open set $V_e$ containing $e$ such that $V_e^{-1} \subseteq U$. Let $V = \bigcup\{V_e : e \in E(S)\}$. Clearly, $V^{-1} \subseteq U$. Let $x \in \overline{A}$. Then $(xV)^*$ is an open set containing $x$ in $S$. So $(xV)^* \cap A \neq \emptyset$. Let $a \in (xV)^* \cap A$. Then $x^{-1}a \in V$ and $x^0 = a^0$. Let $x^{-1}a = v$, for some $v \in V$. Then $a^{-1}x = v^{-1} \in V^{-1}$ and $a^0 = x^0$. It follows that $x \in (aV^{-1})^* \subseteq (AV^{-1})^* \subseteq (AU)^*$. Therefore, $\overline{A} \subseteq (AU)^*$.
\end{proof}

\begin{theorem}
Let $S$ be a band of topological groups and $x \in S$. If $\mathcal{B}$ is a base at $x^0$ in $S$, then $\mathcal{U} = \{(xU)^* : U \in \mathcal{B}\}$ is a base at $x$ in $S$.
\end{theorem}

\begin{proof}
Let $W$ be an open set containing $x$ in $S$. Then $(x^{-1}W)^*$ is an open set containing $x^0$ in $S$. Since $\mathcal{B}$ is a base at $x^0$, there exists an open set $U$ in $\mathcal{B}$ such that $x^0 \in U \subseteq (x^{-1}W)^*$. Then $(xU)^*$ is an open set containing $x$ in $S$. Now, we show that $(xU)^* \subseteq W$. Let $u \in (xU)^*$. Then $x^{-1}u \in U$ and $x^0 = u^0$. It follows that $x^{-1}u \in (x^{-1}W)^*$ and this implies that $u \in W$. Hence $x \in (xU)^* \subseteq W$, where $(xU)^* \in \mathcal{U}$. Therefore, $\mathcal{U}$ is a base at $x$ in $S$.
\end{proof}

\begin{cor}
Let $S$ be a band of topological groups and $x \in S$. If $\mathcal{B}$ is a base at $x^0$ in $S$, then $\mathcal{U} = \{(Ux)^* : U \in \mathcal{B}\}$ is a base at $x$ in $S$.
\end{cor}

\begin{theorem} \label{th 3.6}
Let $S$ be a band of topological groups and for each $e \in E(S)$, let $\mathcal{U}_e$ be a base at $e$ in $S$. Then

(i) for each open set $U$ and for every $x\in U$, there exists $V \in \mathcal{U}_{x^0}$ such that $(Vx)^* \subseteq U$.

(ii) for every $W \in \mathcal{U}_e$, if $(xy)^0 \in W$, there exist $U \in \mathcal{U}_{x^0}$ and $V \in \mathcal{U}_{y^0}$ such that $(UV)^0 \subseteq W$,

(iii) for every $U \in \mathcal{U}_e$, there is an element $V \in \mathcal{U}_e$ such that $V^2 \subseteq U$,

(iv) for every $U \in \mathcal{U}_e$ and every $y \in U \cap H_{_e}$, there is an element $V \in \mathcal{U}_e$ such that $(yVy^{-1})^* \subseteq U$,

(v) $\{e\} = \cap \mathcal{U}_e$, if $S$ is Hausdorff.
\end{theorem}

\begin{proof}
$(i)$ Let $U$ be an open set and $x \in U$. Then $x^0x \in U$. Now, the mapping $(x,y) \mapsto xy$ being continuous, there exists $V\in \mathcal{U}_{x^0}$ such that $Vx \subseteq U$ and thus $(Vx)^* \subseteq Vx \subseteq 
U$.

$(ii)$ Let $W \in \mathcal{U}_e$ and $(xy)^0 \in W$. Then $(x^0y^0)^0\in W$. Now, the mapping $x \mapsto x^0$ being continuous, there exists an open set $W_1$ containing $x^0y^0$ such that $(W_1)^0 \subseteq W$. Since $x^0y^0 \in W_1$ and the mapping $(x,y) \mapsto xy$ is continuous, there exist $U \in \mathcal{U}_{x^0}$ and $V \in \mathcal{U}_{y^0}$ such that $UV \subseteq W_1$ and hence $(UV)^0 \subseteq W$. 

$(iii)$ Follows from the fact that the mapping $(x,y) \mapsto xy$ is continuous.

$(iv)$ Let $U \in \mathcal{U}_e$ and $y \in U \cap H_{_e}$. Then $yey^{-1} = e \in U$ and hence there exists $V \in \mathcal{U}_e$ such that $yVy^{-1} \subseteq U$ and therefore, clearly $(yVy^{-1})^* \subseteq yVy^{-1} \subseteq U$. 

$(v)$ Let $S$ be Hausdorff. Clearly, $\{e\} \subseteq \cap \mathcal{U}_e$. Let $y \in \cap \mathcal{U}_e$. If possible, let $e \neq y$. Then there exist two disjoint open sets $U$ and $V$  in $S$ containing $e$ and $y$ respectively. Since $\mathcal{U}_e$ is a base at $e$, there is an element $W \in \mathcal{U}_e$ such that $e \in W \subseteq U$. Then $y \in W$ and thus $U \cap V \neq \emptyset$, a contradiction. Therefore, $y = e$ and hence $\{e\} = \cap \mathcal{U}_e$.
\end{proof}

\begin{theorem}
Let $S$ be a band of topological groups and for each $e \in E(S)$, let $\mathcal{U}_e$ be a base at $e$ in $S$. Then

(1) for every $U \in \mathcal{U}_e$, there is an element $V \in \mathcal{U}_e$ such that $V^{-1} \subseteq U$;

(2) for every $U \in \mathcal{U}_e$ and every $y \in U$, there is an element $V \in \mathcal{U}_{y^0}$ such that $(Vy)^* \subseteq U$;

(3) for every $U \in \mathcal{U}_e$ and every $x,y \in S$ with $x^{-1}y \in U$, there exists $V \in \mathcal{U}_{y^0}$ such that $x^{-1}Vy \subseteq U$;

(4) for every $a, b \in S$ and $W\in \mathcal{U}_{(ab)^0}$, there exist $U \in \mathcal{U}_{a^0}$ and $V \in \mathcal{U}_{b^0}$ such that $UaVb(ab)^{-1} \subseteq W$.

(5) for $U,V \in \mathcal{}{U}_e$, there is $W \in \mathcal{U}_e$ such that $W \subseteq U \cap V$.

Conversely, let $S$ be a cryptogroup and for each $e \in E(S)$, let $\mathcal{U}_e$ be a family of subsets of $S$ satisfying conditions (1) - (5). Let $\mathcal{U} = \bigcup\limits_{e \in E(S)}\mathcal{U}_e$. Then the family $\mathcal{B}_{\mathcal{U}} = \{(Ua)^* : U \in \mathcal{U} \, \, \mbox{and} \, \, a\in S\}$ is a base for some topology $\tau$ on $S$ and $(S, \tau)$ is a band of topological groups.
\end{theorem}

\begin{proof}Let $S$ be a band of topological groups.

{\it (1)} Follows from the fact that the mapping  $x \mapsto x^{-1}$ is continuous.

{\it (2)} Follows from Theorem \ref{th 3.6} (i).

{\it (3)} Let $U \in \mathcal{U}_e$ and $x,y \in S$ with $x^{-1}y \in U$. Then $x^{-1}y^0y \in U$. Now, the mapping $(x,y) \longrightarrow xy$ being continuous, there exists $V\in \mathcal{U}_{y^0}$ such that $x^{-1}Vy \subseteq U$.

{\it (4)} Let $a, b \in S$ and $W\in \mathcal{U}_{(ab)^0}$. Then $ab(ab)^{-1} \in W$, i.e., $a^0ab^0b(ab)^{-1} \in W$. Now, the mapping $(x,y) \longrightarrow xy$ being continuous, there exist $U \in \mathcal{U}_{a^0}$ and $V \in \mathcal{U}_{b^0}$ such that $UaVb(ab)^{-1} \subseteq W$.

{\it (5)} Follows from the fact that $\mathcal{U}_e$ is a base at $e$ of $S$.

\vspace{1em}
For the converse part, let $S$ be a cryptogroup and for each $e \in E(S)$, let $\mathcal{U}_e$ be a family of subsets of $S$ satisfying conditions {\it (1) - (5)} in this theorem. Let $\mathcal{U} = \bigcup\limits_{e \in E(S)}\mathcal{U}_e$ and let $\tau$ be the collection of all subsets $W$ of $S$ such that for each $x \in W$, there esists $U \in \mathcal{U}_{x^0}$ such that $(Ux)^* \subseteq W$.

Step 1: $\tau$ is a topology on $S$.

Clearly, $\emptyset, S \in \tau$ and for any subfamily $\gamma$ of $\tau$, $\bigcup\limits_{W \in \gamma}W \in \tau$. Let $W_1, W_2 \in \tau$. Set $W = W_1 \cap W_2$. Let $x \in W$. Then $x \in W_1$ and $x \in W_2$. This implies that there exist $U_1, U_2 \in \mathcal{U}_{_{x^0}}$ such that $(U_1x)^* \subseteq W_1$ and $(U_2x)^* \subseteq W_2$. By $(5)$, it follows that  there is an element $U\in \mathcal{U}_{x^0}$ such that $U \subseteq  U_1 \cap U_2$. This implies that $(Ux)^* \subseteq (U_1x)^* \cap (U_2x)^* \subseteq W_1 \cap W_2 = W$. Therefore, $W \in \tau$ and hence $\tau$ is a topology on $S$.

Step 2: For any $x\in S$ and each $U \in \mathcal{U}$, $(Ux)^* \in \tau$.

Let $W = (Ux)^*$ and let $y \in (Ux)^*$. Then $yx^{-1} \in U$, $y^0 = x^0$ and hence by condition {\it (2)}, there exist $U_1 \in \mathcal{U}_{y^0}$ such that $(U_1yx^{-1})^* \subseteq U$. Now, we show that $(U_1)^* \subseteq W$. Let $z \in (U_1y)^*$. Then $zy^{-1} \in U_1$ and $z^0 = y^0$. Therefore, $zy^{-1} = zx^{-1}xy^{-1}= zx^{-1}(yx^{-1})^{-1} \in U_1$ and $(zx^{-1})^0 =(yx^{-1})^0$. This implies $zx^{-1} \in (U_1yx^{-1})^* \subseteq U$ and $z^0 = x^0$. Therefore, $z \in (Ux)^*$ and hence $(U_1y)^* \subseteq (Ux)^* $.

Step 3: From the construction of $\tau$, it at once follow that $\mathcal{B}_{\mathcal{U}} = \{(Ua)^* : a \in S \, \, \text{and}\,\, U \in \mathcal{U}\}$ forms a base for the topology $\tau$.

Step 4: The mapping $\mu : S\times S \longrightarrow S$ defined by $\mu (x, y) = xy$, is continuous with respect to the topology $\tau$.

 Let $W$ be an open set in $(S, \tau)$ and $a,b \in S$ with $ab \in W$. Then there exists $W_1 \in \mathcal{U}_{(ab)^0}$ such that $(W_1ab)^* \subseteq W $. Since $(ab)^0 \in W_1$, so by the condition {\it (4)}, there exist $U \in \mathcal{U}_{a^0}$ and $V\in \mathcal{V}_{b^0}$ such that $UaVb(ab)^{-1}\subseteq W_1$. Then $(Ua)^*$ and $(Vb)^*$ are open sets containing $a$ and $b$ respectively. Now, we show that $(Ua)^*(Vb)^* \subseteq W$. For this, let $x \in (Ua)^*$ and $y \in (Vb)^*$. Then $xa^{-1} \in U$, $x^0 = a^0$ and $yb^{-1} \in V$, $y^0 = b^0$. Let $xa^{-1}= u$ and $yb^{-1} = v$, for some $u \in U$ and $v \in V$. Then $x = ua$ and $y = vb$ and thus $xy(ab)^{-1} = uavb(ab)^{-1} \in UaVb(ab)^{-1} \subseteq W_1$ together with $(ab)^0 = (a^0b^0)^0 = (x^0y^0)^0 = (xy)^0$ imply that $xy \in (W_1ab)^* \subseteq W$. Hence $(Ua)^*(Vb)^* \subseteq W$ and therefore, $(S, \tau)$ is a topological semigroup.
 
Step 5: For each $U \in \mathcal{U}$ and $x \in S$, $(xU)^* \in \tau$.
 
Let $U \in \mathcal{U}$ and $x \in S$. Set $W = (xU)^*$. Let $y \in W$. Then $x^{-1}y \in U$ and $x^0 = y^0$. By condition {\it (3)}, there exists $V \in \mathcal{U}_{y^0}$ such that $x^{-1}Vy \subseteq U$. We now show that $(Vy)^* \subseteq W$. Let $z \in (Vy)^*$. Then $zy^{-1} \in V$ and $z^0 = y^0$. Let $zy^{-1} = v$, where $v \in V$. Then $z = vy$ and thus $x^{-1}z = x^{-1}vy \in x^{-1}Vy \subseteq U$ such that $x^0 = z^0$. Therefore, it follows that $z \in (xU)^*$. Therefore, $(Vy)^* \subseteq (xU)^*$ and hence $(xU)^* \in \tau$.
 
Step 6: The mapping $\gamma : S \longrightarrow S$ defined by $\gamma (x) = x^{-1}$, is continuous with respect to the topology $\tau$.
 
Let $W$ be an open set in $(S,\tau)$ and $x^{-1} \in W$ for some $x \in S$. Then by condition {\it (2)}, there exists $W_1 \in \mathcal{U}_{x^0}$ such that $(W_1x^{-1})^* \subseteq W$. Now, since $x^0 \in W_1$, so by condition {\it (1)}, there exists $V \in \mathcal{U}_{x^0}$ such that $V^{-1} \subseteq W_1$. This implies that $(V^{-1}x^{-1})^* \subseteq (W_1x^{-1})^* \subseteq W$. Now, we show that $(V^{-1}x^{-1})^* = ((xV)^*)^{-1}$. Now, $z\in (V^{-1}x^{-1})^*$ if and only if $zx \in V^{-1}$ such that $z^0 = x^0$ if and only if $x^{-1}z^{-1} \in V$ such that $(z^{-1})^0 = x^0$ if and only if $z^{-1} \in (xV)^*$ if and only if $z \in ((xV)^*)^{-1}$. Thus, we have $(V^{-1}x^{-1})^* = ((xV)^*)^{-1}$ and therefore, $((xV)^*)^{-1} \subseteq W$. In this connection, it is important to note that $(xV)^*$ is an open set containing $x$. Hence $\gamma$ is continuous. Consequently, $(S, \tau)$ is a topological cryptogroup.
 
Step 7: Each $\mathcal{H}$-class is open in $(S, \tau)$.
 
Let $H$ be a $\mathcal{H}$-class and $x \in H$. Let $U \in \mathcal{U}_{x^0}$. Then $x \in (Ux)^* \subseteq H$ and hence $H$ is open in $(S,\tau)$.
 
Therefore, by Corollary \ref{main theorem 2}, it follows that $(S, \tau)$ is a band of topological groups. 
\end{proof}

\begin{theorem}
Let $f: S \longrightarrow T$ be a homomorphism of bands of topological groups. If $f$ is continuous at each idempotent element of $S$, then $f$ is continuous.
\end{theorem}

\begin{proof}
Let $W$ be an open subset of $T$ and let $x \in f^{-1}(W)$. Then $f(x) \in W$, i.e., $f(x)(f(x))^0 \in W$. Now, $T$ being a band of topological groups, there exists an open set $V$ containing $(f(x))^0$ in $T$ such that $f(x)V \subseteq W$ and thus $(f(x)V)^* \subseteq f(x)V$. Since $(f(x))^0 = f(x^0)$, it follows that $V$ is an open set containing $f(x^0)$ in $T$. Now, $f$ being continuous at $x^0$, there exists an open set $U$ containing $x^0$ such that $f(U) \subseteq V$ and $(xU)^*$ is an open set containing $x$ in $S$. Now, we show that $(xU)^* \subseteq f^{-1}(W)$. For this, let $y \in (xU)^*$. Then $x^{-1}y \in U$ and $x^0 = y^0$. Also, $f(x^{-1}y) \in f(U) \subseteq V$ implies $(f(x))^{-1}f(y) \in V$ and $(f(x))^0 = f(x^0) = f(y^0) = (f(y))^0$. This implies $f(y) \in (f(x)V)^*$ and thus we have $f(y) \in f(x)V \subseteq W$ and so $y \in f^{-1}(W)$. Hence $(xU)^* \subseteq f^{-1}(W)$ and consequently, $f$ is continuous.
\end{proof}

\begin{theorem}
Let $S$ be a Hausdorff topological semigroup. Then for any $t \in S$, the set $S_{_t} = \{x \in S : xt = tx\} $ is a closed subsemigroup of $S$.
\end{theorem}

\begin{proof}
Clearly, the set $S_{_t} = \{x \in S : xt = tx\} $ is a subsemigroup of $S$. Let $x \in \overline{S_{_t}}$. If possible, let $x \notin S_{_t}$, i.e., $xt \neq tx$. Then there exist two open sets $W_1$ and $W_2$ in $S$ containing $xt$ and $tx$ respectively such that $W_1\cap W_2 = \emptyset$. Now, $S$ being a topological semigroup, there exists an open set $U$ containing $x$ in $S$ such that $Ut \subseteq W_1$ and $tU \subseteq W_2$. Since $x \in \overline{S_{_t}}$, it follows that $U \cap S_{_t} \neq \emptyset$. Let $u \in U \cap S_{_t}$. Then $tu = ut$  and this implies $W_1 \cap W_2 \neq \emptyset$, a contradiction. Hence $x \in S_{_t}$ and therefore $S_{_t}$ is a closed subsemigroup of $S$.
\end{proof}

We now register an elementary result that follows from \cite{koch}.

\begin{theorem}
Let $S$ be a Hausdorff topological semigroup. Then $E(S)$ is closed in $S$.
\end{theorem}

\begin{theorem}
Let $S$ be a Hausdorff topological semigroup and $n$ be a positive integer. Then the set $S[n] = \{x \in S : x^n \in E(S)\}$ is closed in $S$.
\end{theorem}

\begin{proof}
Let $x \in \overline{S[n]}$. If possible, let $x \notin S[n]$. Then $x^n \notin E(S)$. Since $E(S)$ is closed, there exists an open set $W$ containing $x^n$ such that $W \cap E(S) = \emptyset$. Now, $S$ being a topological semigroup, there exists an open set $U$ containing $x$ such that $U^n \subseteq W$. Since $x \in \overline{S[n]}$, it follows that $S[n] \cap U \neq \emptyset$. Let $y \in S[n] \cap U$. Then $y^n \in U^n \cap E(S)$ and thus $W \cap E(S) \neq \emptyset$, a contradiction. Therefore, $x \in S[n]$ and hence $S[n]$ is closed in $S$.
\end{proof}

\section{Topological properties of a band of topological groups}

In \cite{mp}, we studied the topological properties of semilattice of topological groups. In this section, we study topological properties, like $T_i$-axioms, regularity, complete regularity, normality, locally connectedness, first countability, second countability, separability, metrizability, etc. of a band of topological groups. 

Similar to the proof of \cite[Theorem 3.1]{mp}, one can easily establish the following result.

\begin{theorem} \label{th:3.1}
Let $(S, \tau) = (B; G_{_{\alpha}}, \tau_{_\alpha})$ be a band of topological groups. Then

(i) $(S,\tau)$ is $T_0$ if and only if each $(G_{_{\alpha}}, \tau_{_\alpha})$ is $T_0$.

(ii) $(S, \tau)$ is $T_1$ if and only if each  $(G_{_{\alpha}}, \tau_{_\alpha})$ is $T_1$.

(iii) $(S, \tau)$ is Hausdorff if and only if each  $(G_{_{\alpha}}, \tau_{_\alpha})$ is Hausdorff.

(iv) $(S, \tau)$ is regular (topological property) if and only if each  $(G_{_{\alpha}}, \tau_{_\alpha})$ is regular (topological property).

(v) $(S, \tau)$ is completely regular (topological property) if and only if each  $(G_{_{\alpha}}, \tau_{_\alpha})$ is completely regular (topological property).

(vi) $(S, \tau)$ is normal (topological property) if and only if each  $(G_{_{\alpha}}, \tau_{_\alpha})$ is normal (topological property).
\end{theorem}

Using Theorem \ref{th:3.1}, we at once have the following corollary.

\begin{theorem}\label{cor 3.2}
Let $S$ be a band of topological groups. Then the following conditions are equivalent:

(i) $S$ is $T_0$;

(ii) $S$ is $T_1$;

(iii) $S$ is Hausdorff;

(iv) $S$ is regular (topological property);

(v) $S$ is completely regular (topological property).
\end{theorem}

\begin{remark}
In Corollary \ref{cor 3.2}, we have shown that a Hausdorff band of topological groups is always a completely regular (topological property) band of topological groups. But, in general, a Hausdorff band of topological groups may not be topologically normal. This follows from the following example.
\end{remark}

\begin{example}\cite{th} 
Let $G= \mathbf{Z}^A$, where $\mathbf{Z}$ is the additive Abelian group of all integers endowed with the usual topology and $A$ is an uncountable set. Then the set $G$ endowed with the product topology is a Hausdorff topological group which is not normal. 
\end{example} 

Before going to further study, we first state two theorems.

\begin{theorem}\label{th: locally compact normal}\cite{th}
Every locally compact Hausdorff topological group is normal.
\end{theorem}

Since each $\mathcal{H}$-class of a band of topological groups is closed, applying Theorem \ref{th:3.1} and Theorem \ref{th: locally compact normal}, we have the following result.

\begin{theorem}\label{th:loc con semigr}
Every locally compact Hausdorff band of topological groups is normal.
\end{theorem}

Similar to the proof of  \cite[Theorem 3.9]{mp}, we can prove the following result.

\begin{theorem} 
Let $S$ be a band of topological groups. Then $S$ is locally connected if and only if each $\mathcal{H}$-class is locally connected.
\end{theorem}

Now, we recall a result that follows from \cite{th}

\begin{theorem} \label{th:loc con}
Let $G$ be a locally compact Hausdorff topological group. Also, assume that $C$ is a closed subset of $G$ and $U$ is an open set containing $C$. Then there is a continuous function $f : G \longrightarrow [0, 1]$ such that $f(C) = \{1\}$ and $f(G \setminus U) = \{0\}$
\end{theorem}

Combining Theorem \ref{th:loc con semigr} and Theorem \ref{th:loc con}, we have the following corollary.

\begin{cor}
Let $S$ be a locally compact Hausdorff band of topological groups. For any closed set $C$ in $S$ and any open set $U$ containing $C$, there is a continuous function $f : S \longrightarrow [0,1]$ such that $f(C) = \{1\}$ and $f(S \setminus U) = \{0\}$.
\end{cor}

\begin{theorem}
Let $S$ be a band of topological groups in which $E(S)$ is countable. Then $S$ is separable if and only if each $\mathcal{H}$- class is separable.  
\end{theorem}

\begin{proof}
Let $S$ be separable and $K$ be a countable dense subset of $S$. Let $H$ be a $\mathcal{H}$-class of $S$. Let $A = K \cap H$. Clearly, $A$ is a countable subset of $H$. To show $H$ is separable, it suffices to show that $\overline{A} = H$. Now, $A \subseteq H$ implies $\overline{A} \subseteq \overline{H}$. Since $S$ is a band of topological groups, $H$ is closed in $S$ and so $\overline{A} \subseteq H$. For the reverse inclusion, let $h\in H$. Let $U$ be an open set containing $h$ in $S$. Since $S$ is a band of topological groups, $H$ is open in $S$ and hence $U\cap H$ is an open set containing $h$ in $S$. Since $\overline{K} = S$, it follows that $(U \cap H) \cap K \neq \emptyset$, i.e., $U \cap  A \neq \emptyset$. This implies $h \in \overline{A}$. Therefore, $\overline{A} = H$ and consequently, $H$ is separable. 

Conversely, let each $\mathcal{H}$-class be separable. Since $E(S)$ is countable and $E(S)$ is homeomorphic to $S/\mathcal{H}$, $S/\mathcal{H}$ is a countable set. Let $\{H_i : i \geq 1\}$ be the collection of all $\mathcal{H}$-classes of $S$. For each $i \geq 1$, let $A_i$ be a countable dense subset of $H_i$. Let $B = \bigcup\limits_{i \geq 1}A_i$. Then $B$ is a countable subset of $S$. Now, we show that $\overline{B} = S$. Let $s \in S$ and $U$ be an open set containing $s$ in $S$. Then $s \in H_j$, for some $j \geq 1$. Since $\overline{A_j} = H_j$, $U \cap A_j \neq \emptyset$ and this implies that $U \cap B \neq \emptyset$. Therefore, it follows that $s \in \overline{B}$ and hence $\overline{B} = S$. Consequently, $S$ is separable.
\end{proof}

\begin{theorem}\label{first countable}
Let $S$ be a band of topological groups. Then $S$ is first countable if and only if each $\mathcal{H}$-class is first countable.
\end{theorem}

\begin{proof}
If $S$ is first countable, then clearly each $\mathcal{H}$-class is first countable. 

Conversely, let each $\mathcal{H}$-class be first countable. Let $x \in S$. Then $x \in H_{_x}$. Since $H_{_x}$ is first countable, $x$ has a countable local basis, say $\{B_{_i} : i \geq 1\}$ in $H_{_x}$. Again, since $S$ is a band of topological groups, $H_{_x}$ is open in $S$ and so  $\{B_{_i}: i \geq 1\}$ is a countable collection of open sets in $S$. Let $G$ be an open set containing $x$ in $S$. Then $G \cap H_{_x}$ is an open set containing $x$ in $H_{_x}$ and therefore, there exists some $B_j$ such that $x \in B_j \subseteq G \cap H_{_x}$ and this implies $x \in B_{_j} \subseteq G$. This proves that $x$ has a countable local basis $\{B_{_i} : i \geq 1\}$ in $S$ and hence $S$ is first countable.
\end{proof}

\begin{theorem}
Let $S$ be a band of topological groups in which $E(S)$ is countable. Then $S$ is second countable if and only if each $\mathcal{H}$-class is second countable.
\end{theorem}

\begin{proof}
Let $S$ be second countable and $\{B_i : i \geq 1\}$ is the countable base for $S$. Let $H$ be a $\mathcal{H}$-class of $S$. Now $\{B_i \cap H: i \geq 1 \}$ is a countable collection of open sets in $H$. Let $G$ be an open set in $H$. Then $G = V \cap H$, for some open set $V$ in $S$. Since $\{B_i : i \geq 1\}$ is a base for $S$, $B_j \subseteq V$, for some $j \geq 1$. From this, it follows that $B_j \cap H \subseteq V \cap H = G$. Therefore, $\{B_i \cap H: i \geq 1 \}$ is a countable basis of $H$ and hence $H$ is second countable. 

Conversely, let each $\mathcal{H}$-class be second countable. Since $E(S)$ is countable and $E(S)$ is homeomorphic to $S/\mathcal{H}$, $S/\mathcal{H}$ is a countable set. Let $\{H_i : i \geq 1\}$ be the collection of all $\mathcal{H}$-classes of $S$. For each $i \geq 1$, let $\{C_{ij} : j \geq 1\}$ be a countable base for $H_i$. Now, $S$ being a band of topological groups, $H_i$ is open in $S$, for all $i \geq 1$. Therefore, $\{C_{ij}: i,j \geq 1\}$ is a countable collection of open sets in $S$. Let $U$ be an open set in $S$ and $x \in U$. Then $x \in H_k$, for some $k \geq 1$. Now, $U \cap H_k$ is an open set containing $x \in H_k$ and $\{C_{kj} : j\geq 1\}$ is a basis of $H_k$. So $x \in C_{ks} \subseteq U \cap H_k $, for some $s \geq 1$. This implies $x \in C_{ks} \subseteq U$ and therefore, $\{C_{ij} : i,j \geq 1\}$ is a countable basis of $S$. Consequently, $S$ is secondly countable.
\end{proof}

\begin{theorem}\label{metrizable}
A band of topological groups $(S,\tau)$ is metrizable if and only if each $\mathcal{H}$-class is metrizable.
\end{theorem}

\begin{proof}
If $(S, \tau)$ is metrizable, then clearly each $\mathcal{H}$-class is metrizable. 

Conversely, let each $\mathcal{H}$-class be metrizable. Let $S = \bigcup\limits_{\alpha \in B}G_{_\alpha}$, where $B$ is a band and each $G_{_\alpha}$ is a $\mathcal{H}$-class. For each $\alpha \in B$, let $d_{_\alpha}$ be the metric on $G_{_\alpha}$ for which $G_{_\alpha}$ is metrizable. Consider the standard bounded metric $\bar{d_{_\alpha}}$ on $G_{_\alpha}$ defined by $\bar{d_{_\alpha}}(x,y) = \text{min}\{1, \,d(x,y)\}$, for all $\alpha$. For each $\alpha$, let $\tau_{_\alpha}$ be the induced topology on $G_{_\alpha}$. Then $d_{_\alpha}$ and $\bar{d_{_\alpha}}$ induces the same topology $\tau_{_\alpha}$ on $G_{_\alpha}$, for all $\alpha$. Define a metric $d$ on $S$ as follows: $x \in G_{_\alpha}, y \in G_{_\beta}$,
 
\begin{equation}
d(x,y) =
\begin{cases}
\bar{d_{_\alpha}}(x,y) & \text{if $\alpha = \beta$ }\\
1 & \text{otherwise}
\end{cases}       
\end{equation}
Clearly, $d$ is a metric on $S$. Let $\tau_{_d}$ be the topology on $S$ corresponding to the metric $d$ on $S$. We show that $\tau_{_d} = \tau $. Let $x \in S$ and $\epsilon > 0$. Then $x \in G_{_\alpha}$, for some $\alpha \in B$. Now, $B(x, \epsilon) = \{y \in S : d(x,y) < \epsilon\}$. For $\epsilon \leq 1$, $B(x, \epsilon) = \overline{B_{_\alpha}}(x, \epsilon)$ and for $\epsilon >1$, $B(x, \epsilon) = (\bigcup\limits_{\alpha \neq \beta}G_{_\beta}) \bigcup \overline{B_{_\alpha}}(x, \epsilon)$, where $\overline{B_{_\alpha}}(x, \epsilon) = \{y \in S : \bar{d_{_\alpha}}(x,y) < \epsilon\}$. Now, $\overline{B_{_\alpha}}(x, \epsilon)$ is open in $G_{_\alpha}$ and hence open in $S$. Therefore, $B(x, \epsilon)$ is open in $(S, \tau)$ and thus $\tau_{_d} \subseteq \tau$. For the reverse inclusion, let $G$ be an open set in $(S, \tau)$ and $y \in G$. Then $y \in G_{_\beta}$, for some $\beta \in B$ and so $G \cap G_{_\beta} \in \tau_{_\beta}$. Now, $G_{_\beta} = \bigcup\limits_{n\geq 1}B(y, 1/n) \in \tau_{_d}$ and therefore, it follows that $G \cap G_{_\beta} \in \tau_{_d}$. Therefore, there is a $\delta_{_y} >0$ such that $ y \in B(y, \delta_{_y}) \subseteq G \cap G_{_\beta}$. This implies $G = \bigcup\limits_{y \in G} B(y, \delta_{_y})$ and so $G \in \tau_{_d}$. Therefore, $\tau = \tau_{_d}$ and consequently, $S$ is metrizable. 
\end{proof}

The following well-known result can be found, for instance, in \cite{bir}.

\begin{theorem}\label{groupmetrizability}
A Hausdorff topological group G is metrizable if and only if it satisfies the first countability axion.
\end{theorem}

\begin{theorem}
A Hausdorff band of topological groups $S$ is metrizable if and only if it is first countable.
\end{theorem}

\begin{proof}
The necessary part is obvious. For the converse part, let $S$ be first countable. Then each $\mathcal{H}$-class is first countable. Since $S$ is a Hausdorff band of topological groups, it at once follows that each 
$\mathcal{H}$-class is a Hausdorff topological group. Therefore, by Theorem \ref{groupmetrizability}, each $\mathcal{H}$-class is metrizable and by Theorem \ref{metrizable}, it follows that $S$ is metrizable.
\end{proof}

\section{Full subcryptogroup}

In this section, we discuss the topological properties of full subcryptogroups of a band of topological groups. Before going to study the properties of a full subcryptogroup, we concentrate on a symmetric subset of a band of topological groups.

\begin{theorem}
Let $S$ be a band of topological groups. Then for each $e \in E(S)$, $S$ has an open base at $e$ of symmetric neighborhoods.
\end{theorem}

\begin{proof}
Let $e \in E(S)$ and $U$ be an arbitrary open set containing $e$ in $S$. Set $V = U \cap U^{-1}$. Then $V$ is an open set containing $e$ with $V^{-1} = V$. Moreover, $V \subseteq U$. Hence the result.
\end{proof}

\begin{theorem}
Let $S$ be a cryptogroup endowed with a topology $\tau$ such that the inverse mapping $x \mapsto x^{-1}$ is continuous. Then for any symmetric subset $A$ of $S$, $\overline{A}$ is also symmetric.
\end{theorem}

\begin{proof}
Let $A$ be a symmetric subset of $S$. We show that $(\overline{A})^{-1} = \overline{A}$. Let $x \in (\overline{A})^{-1}$. Then $x^{-1} \in \overline{A}$. Let $U$ be an open set containing $x$ in $S$. Since the mapping $x \mapsto x^{-1}$ is continuous, it is a homeomorphism. This implies that $U^{-1}$ is an open set containing $x^{-1}$ in $S$. Since $x^{-1} \in \overline{A}$, $U^{-1} \cap A \neq \emptyset$, it follows that $U \cap A^{-1} \neq \emptyset$, i.e., $U \cap A \neq \emptyset$. This implies $x \in \overline{A}$ and thus $(\overline{A})^{-1} \subseteq \overline{A}$. For the reverse inclusion, let $y \in \overline{A}$ and let $V$ be an open set containing $y^{-1}$ in $S$. Then $V^{-1}$ is an open set containing $y$ in $S$. This implies that $V^{-1} \cap A \neq \emptyset$ and thus $V \cap A^{-1} \neq \emptyset$, i.e., $V \cap A \neq \emptyset$. This implies $y^{-1} \in \overline{A}$, i.e., $y \in (\overline{A})^{-1}$. Therefore, $\overline{A} \subseteq (\overline{A})^{-1}$  and consequently, $(\overline{A})^{-1} = \overline{A}$.
\end{proof}

Recall that the minimal cardinality \cite{Ar} of a local base at a point $x$ of a topological space $X$ is called the character of $X$ at $x$ and is denoted by $\chi(x, X)$.

\vspace{1em}
The following well-known result can be found, for instance, in \cite{Ar}.

\begin{theorem}\label{dense countable}
If $Y$ is a dense subspace of a regular topological space $X$, then $\chi (y, Y) = \chi (y, X)$, for all $y \in Y$.
\end{theorem}

\begin{lemma}\label{cryptosubgroup}
Closure of a subcryptogroup of a topological cryptogroup $S$ is also a subcryptogroup of $S$.
\end{lemma}

\begin{proof}
Let $N$ be a subcryptogroup of a topological cryptogroup $S$. Let $x, y \in \overline{N}$. If possible, let $xy \notin \overline{N}$. Then there exists an open set $W$ containing $xy$ in $S$ such that $W \cap N = \emptyset$. Now $S$ being a topological cryptogroup, there exist open sets $U$ and $V$ containing $x$ and $y$ respectively in $S$ such that $UV \subseteq W$. Now, $x,y \in \overline{N}$ implies that $U \cap N \neq \emptyset$ and $V \cap N \neq \emptyset$. Let $p \in U \cap N$ and $q \in V \cap N$. Then $pq \in N^2 \subseteq N$. Again $pq \in UV \subseteq W$. This implies that $W \cap N \neq \emptyset$, a contradiction. So $xy \in \overline{N}$. Now, we show that $x^{-1} \in \overline{N}$. If possible, let $x^{-1} \notin \overline{N}$. Then there is an open set $G$ containing $x^{-1}$ in $S$ such that $G \cap N = \emptyset$. This implies $G^{-1}$ is an open set containing $x$. Since $x \in \overline{N}$, $N \cap G^{-1} \neq \emptyset$ and this implies $N^{-1} \cap G \neq \emptyset$, i.e., $N \cap G \neq \emptyset$, a contradiction. Therefore, $x^{-1} \in \overline{N}$ and thus $\overline{N}$ is a subcryptogroup of $S$. 
\end{proof}

\begin{theorem}
Let $S$ be a Hausdorff band of topological groups and $H$ be a first countable full subcryptogroup of $S$. Then $\overline{H}$ is also a first countable full subcryptogroup of $S$.
\end{theorem}

\begin{proof}
Clearly, $\overline{H}$ is a full subcryptogroup of $S$. Since $S$ is a Hausdorff band of topological groups, $S$ is regular (topological property). Let $x \in \overline{H}$. Then $x^0 \in E(S) \subseteq H$.  Since $H$ is first countable, $\chi(x^0 , H) \leq \aleph_0$ and so by Theorem \ref{dense countable}, it follows that $\chi(x^0, \overline{H}) \leq \aleph_0$. We know that if $\mathcal{B}$ is a local base for $\overline{H}$ at $x^0$, then $\{(Vx)^* : V \in \mathcal{B}\}$ is a local base for $\overline{H}$ at $x$ and therefore, $\chi (x, \overline{H}) \leq \aleph_0$. Hence the theorem.
\end{proof}

 Similar to the proof of \cite[Theorem 7.10]{mp}, we can prove the following result.

\begin{theorem}\label{open implies closed}
Every open full subcryptogroup $K$ of a band of topological groups $S$ is closed in $S$.
\end{theorem}

\begin{theorem}\label{discrete}
If $K$ is a discrete full subcryptogroup of a Hausdorff band of topological groups $S$, then $K$ is closed in $S$.
\end{theorem}

\begin{proof}
Clearly, $\overline{K}$ is a full subcryptogroup of $S$.  Since $K$ is discrete in itself and $S$ is Hausdorff, it follows that $K$ is open in $\overline{K}$. Now, $\overline{K}$ is a topological cryptogroup with respect to the induced topology. Let $a\in \overline{K}$. Also, let $H_a$, $(H_a)_{\overline{K}}$ respectively be the $\mathcal{H}$-class containing the element $a$ in $S$, $\overline{K}$. We now show that $(H_a)_{\overline{K}} = H_a \cap {\overline{K}}$. Clearly, $(H_a)_{\overline{K}} \subseteq  H_a \cap {\overline{K}}$. For the reverse inclusion, let $b \in H_a \cap {\overline{K}}$. Then $a^0 = b^0$. This implies that $a = (ab^{-1})b, b = (ba^{-1})a$ and $a = b(b^{-1}a), b = a(a^{-1}b)$. Since $\overline{K}$ is a subcryptogroup, $a,b \in {\overline{K}}$ implies $ab^{-1}, b^{-1}a \in {\overline{K}}$ and therefore, $b \in (H_a)_{\overline{K}}$. Hence $(H_a)_{\overline{K}} = H_a \cap {\overline{K}}$. Since $S$ is a band of topological groups, $H_a$ is open in $S$ and so $(H_a)_{\overline{K}}$ is open in $\overline{K}$. Therefore, $\overline{K}$ is a band of topological groups and so by Theorem \ref{open implies closed}, it follows that $K$ is closed in $\overline{K}$ and hence $K$ is closed in $S$.
\end{proof}

\begin{theorem}
If $K$ is a locally compact  full subcryptogroup of a Hausdorff band of topological groups $S$, then $K$ is closed in $S$.
\end{theorem}

\begin{proof}
Clearly, $\overline{K}$ is a full subcryptogroup of $S$. Now $K$ being a dense locally compact subspace of $\overline{K}$, by \cite[Theorem 3.3.9]{en}, it follows that $K$ is open in $\overline{K}$. Since $\overline{K}$ is a band of topological groups, so by Theorem \ref{open implies closed}, it follows that $K$ is closed in $\overline{K}$. Therefore, $\overline{K} = K$ and hence $K$ is closed in $S$.
\end{proof}

\begin{theorem}\label{countable compact}
Every discrete full subcryptogroup $H$ of a Hausdorff countably compact band of topological groups $S$ is finite.
\end{theorem}

\begin{proof}
By Theorem \ref{discrete}, we have $H$ is closed in $S$. Since $S$ is countably compact and $H$ is closed in $S$, it follows that $H$ is countably compact. Since $H$ is discrete, we must have $H$ is finite. 
\end{proof}

\begin{theorem}
Every discrete full subcryptogroup $H$ of a Hausdorff Lindel$\ddot{o}$f band of topological groups is countable.
\end{theorem}

\begin{proof}
Follows from the similar argument as in Theorem \ref{countable compact}.
\end{proof}

\begin{definition}
Let $S$ be a band of topological groups and $U$ be an open set in $S$. Then a subset $A$ of $S$ is said to be $U$-disjoint if $q \notin (pU)^*$, for any two distinct elements $p, q \in A$.
\end{definition}

\begin{theorem}\label{U-disjoint}
Let $S$ be a band of topological groups and $E(S)$ be a compact subsemigroup of $S$.
Let $U$ and $V$ be open sets containing $E(S)$ in $S$ such that $V^4 \subseteq U$ and $V^{-1} = V$. If a subset $A$ of $S$ is $U$-disjoint, then the family of open sets $\{(aV)^* : a \in A\}$ is discrete in $S$. 
\end{theorem}

\begin{proof}
Let $x \in S$. Then $(xV)^*$ is an open set containing $x$ in $S$. We show that $(xV)^*$ intersects at most one element of the family $\{(aV)^*: a \in A\}$. If possible, let $(xV)^*$ intersect $(aV)^*$ and $(bV)^*$, for some two distinct elements $a, b$ in $A$. Let $p \in (xV)^* \cap (aV)^*$ and $q \in (xV)^* \cap (bV)^*$. Then $x^{-1}p, a^{-1}p, x^{-1}q, b^{-1}q \in V$ and $x^0 = a^0 = p^0 =q^0 = b^0$. From this it follows that $a^{-1}b = (a^{-1}p)(p^{-1}x)(x^{-1}q)(q^{-1}b) = (a^{-1}p)(x^{-1}p)^{-1}(x^{-1}q)(b^{-1}q)^{-1} \in VV^{-1}VV^{-1} \subseteq V^4 \subseteq U$. Also, $a^0 = b^0$. This implies that $b\in (aU)^*$, which contradicts the fact that $A$ is $U$-disjoint. Hence $(xV)^*$ intersects at most one element of the family $\{(aV)^* : a \in A\}$ and therefore, the family of open sets $\{(aV)^* : a \in A\}$ is discrete in $S$.
\end{proof}

 Before going to prove the last result in this section, we state an important theorem from \cite{koch}. 

 \begin{theorem} \label{th:Wallace}
 Let $X, Y,$ and $Z$ be topological spaces, $A$ be a compact subset of $X$, $B$ be a compact subset of $Y$, $f: X \times Y \longrightarrow Z$ be a continuous function, and $W$ be an open subset of $Z$ containing $f(A \times B)$. Then there exists an open set $U$ in $X$ and an open set $V$ in $Y$ such that $A \subseteq U$, $B\subseteq V$, and $f(U \times V) \subseteq W$.
 \end{theorem}

\begin{theorem}
Let $S$ be a Hausdorff pseudo-compact band of topological groups in which $E(S)$ is a compact subsemigroup of $S$. Then every discrete full subcryptogroup of $S$ is finite.
\end{theorem}

\begin{proof}
Let $H$ be a discrete full subcryptogroup of $S$ and $U$ be an open set in $S$ such that $E(S) = U \cap H$. Since $E(S)$ is a compact subsemigroup of $S$, then by Theorem \ref{th:Wallace}, there exists an open set $V$ containing $E(S)$ such that $V^4 \subseteq U$ and $V^{-1} = V$. First, we show that $H$ is $U$-disjoint. Let $a,b \in H$  with $a \neq b$. If possible, let $b \in (aU)^*$. Then $a^{-1}b \in U$ and $a^0 = b^0$. Since $H$ is a subcryptogroup of $S$, so $a,b \in H$ implies $a^{-1}b \in H$. This implies $a^{-1}b \in U\cap H = E(S)$. Now, $a^0 = b^0$ and $a^{-1}b \in E(S)$ imply that $a = b$, contradiction. Therefore, $b \notin (aU)^*$, for any two distinct elements $a,b$ of $H$. Hence $H$ is $U$-disjoint. So by Theorem \ref{U-disjoint}, the family $\{(aV)^* : a \in H\}$ is discrete in $S$. Now, since $S$ is a Hausdorff band of topological groups, $S$ is completely regular. Again, since $S$ is pseudocompact, it follows that $\{(aV)^* : a \in H\}$ is finite. We claim that $H$ is finite. If not, then there exist $x,y \in H$ with $x \neq y$ such that $(xV)^* = (yV)^*$. Then $x^{-1} y \in V$ and $x^0 = y^0$. This implies that $x^{-1}y \in V^4 \subseteq U$ with $x^0 = y^0$ and thus $y \in (xU)^*$, which contradicts the fact that $H$ is $U$-disjoint. Therefore, $H$ is finite.
\end{proof}

\section{Full normal subcryptogroup}

In this section, we show that the quotient space of a band of topological groups by a full normal subcryptogroup is again a band of topological groups. For this purpose, let us start with some elementary results.

\begin{lemma}
Closure of a full normal subcryptogroup of a topological cryptogroup $S$ is also a full normal subcryptogroup of $S$.
\end{lemma}

\begin{proof}
Let $N$ be a full normal subcryptogroup of a topological cryptogroup $S$. Then $E(S) \subseteq N \subseteq \overline{N}$. By Lemma \ref{cryptosubgroup}, $\overline{N}$ is a  subcryptogroup of $S$. Let $s \in S$ and $t \in \overline{N}$. If possible, let $sts^{-1} \notin \overline{N}$. Then there is an open set $W$ containing $sts^{-1}$ in $S$ such that $W \cap N = \emptyset$. Since $sts^{-1} \in W$, there exist open sets $U$ and $V$ containing $s$ and $t$ respectively in $S$ such that $UVU^{-1} \subseteq W$. Now, $t \in \overline{N}$ implies that $V\cap N \neq \emptyset$. Let $n \in V\cap N$. Then $sns^{-1} \in W\cap N$, a contradiction. This contradiction ensures that $sts^{-1} \in \overline{N}$ and consequently, $\overline{N}$ is a full normal subcryptogroup of $S$.
\end{proof}

\begin{theorem}\label{quotient crypto}
Let $S$ be a cryptogroup and $N$ be a full normal subcryptogroup of $S$. Define a relation $\rho_{_N}$ on $S$ by : for any $a,b \in S$,
\begin{center}
$a\rho_{_N} b$ if and only if $a^{-1}b \in N$ and $a^0 =b^0$.
\end{center}
Then $\rho_{_N}$ is a congruence on $S$ and $S/\rho_{_N}$ is a cryptogroup.
\end{theorem}

\begin{proof}
Clearly, $\rho_{_N}$ is reflexive. Let $a,b \in S$ with $a\rho_{_N} b$. Then $a^{-1}b \in N$ and $a^0 = b^0$. Then $a \, \mathcal{H} \, b$. Now, $H_{_a}$ being a group, $b^{-1}a = (a^{-1}b)^{-1} \in N$ and thus $b\rho_{_N}a$. So, $\rho_{_N}$ is symmetric. For transitivity, let $x, y, z \in S$ with $x\rho_{_N} y$ and $y\rho_{_N}z$. Then $x^{-1}y, y^{-1}z \in N$ and $x^0 = y^0 = z^0$. Now, $x^{-1}z = x^{-1}z^0z = x^{-1}y^0z = (x^{-1}y)(y^{-1}z) \in N$ implies $x \rho_{_N} z$ and hence $\rho_{_N}$ is transitive. Let $p, q, r \in S$ with $p \, \rho_{_N} \, q$. Then $p^{-1}q \in N$ and $p^0 = q^0$. Now, $p^0 = q^0$ implies that $p \, \mathcal{H} \, q$ and then $pr \, \mathcal{H} \, qr$ which implies that $(pr)^0 = (qr)^0$. Again, $(pr)^{-1}(qr)= (pr)^0r^{-1}(rp)^0p^{-1}(pr)^0(qr)= (pr)^0r^{-1}(rp)^0p^{-1}(qr)^0(qr)=(pr)^0r^{-1}(rp)^0p^{-1}(qr)=(pr)^0(r^{-1}((rp)^0(p^{-1}q))r) \in N$. Therefore, $pr \, \rho_{_N} \, qr$ and thus $\rho_{_N}$ is a right congruence on $S$. Also, $p \, \mathcal{H} \, q$ implies $rp \, \mathcal{H} \, rq$ and thus $(rp)^0 = (rq)^0$. Moreover, $(rp)^{-1}(rq)= (rp)^0p^{-1}(pr)^0r^{-1}(rp)^0(rq)= (rp)^0p^{-1}(pr)^0r^{-1}(rq) = (rp)^0p^{-1}(pr)^0q =  (rp)^0p^{-1}\\(qr)^0q = (rp)^0p^{-1}(qr)(qr)^{-1}q = (rp)^0(p^{-1}q)(r(qr)^{-1}q) \in N$, as $r(qr)^{-1}q \in E(S)\subseteq N$. Therefore, $rp \, \rho_{_N} \, rq$ and thus, $\rho_{_N}$ is a left congruence on $S$. Hence $\rho_{_N}$ is a congruence on $S$ and thus $S/\rho_{_N}$ is a semigroup. Now, $S$ being completely regular, $S/\rho_{_N}$ is also completely regular. We note that for any $a \in S$, $a\rho_{_N} = \{b \in S : a^{-1}b \in N\,\, \text{and}\,\, a^0 = b^0\} = \{b \in S : ab^{-1} \in N \,\, \text{and}\,\, a^0 =b^0\}$. We denote the Green's relations $\mathcal{L}, \mathcal{R}$ and $\mathcal{H}$ on $S/ \rho_{_N}$ by $\mathcal{L}_{_N}, \mathcal{R}_{_N}$ and $\mathcal{H}_{_N}$ respectively.  Let $u\rho_{_N}, v\rho_{_N}, w\rho_{_N} \in S/\rho_{_N}$ with $u \rho_{_N} \, \mathcal{H}_{_N} \, v \rho_{_N}$. Then $u \rho_{_N} \, \mathcal{L}_{_N} \, v \rho_{_N}$ and $u \rho_{_N} \, \mathcal{R}_{_N} \, v \rho_{_N}$. Now, $ u \rho_{_N} \, \mathcal{L}_{_N} \, v \rho_{_N}$ implies $u \rho_{_N} = x\rho_{_N} v\rho_{_N}$ and $v \rho_{_N} = y\rho_{_N} u\rho_{_N}$, for some $x\rho_{_N}, y\rho_{_N} \in S/\rho_{_N}$. Then $u (xv)^{-1} = n$ and $v (yu)^{-1}= m$ with $u^0 = (xv)^0$ and $v^0 = (yu)^0$, for some $m, n \in N$. This implies $u = (nx)v$ and $v = (my)u$, i.e.,  $u \, \mathcal{L} \, v$. Similarly, one can show that $u \, \mathcal{R} \, v$. Therefore, $u \, \mathcal{H} \, v$. Since $\mathcal{H}$ is a congruence on $S$, we have $wu \mathcal{H} wv$ and $uw \mathcal{H} vw$. This implies that $w \rho_{_N} u\rho_{_N} \, \mathcal{H}_{_N} \, w \rho_{_N} v\rho_{_N}$ and $u \rho_{_N} w\rho_{_N} \, \mathcal{H}_{_N} \, v \rho_{_N} w\rho_{_N}$. Therefore, $\mathcal{H}_{_N}$ is a congruence on $S/\rho_{_N}$ and consequently, $S/ \rho_{_N}$ is a cryptogroup.
\end{proof}

\begin{notation}
Let $x \rho_{_N}$ be a $\rho_{_N}$-class in $S$ containing an element $x$ in a cryptogroup $S$. Then $x \rho_{_N} = (xN)^*$, where $(xN)^* = \{y \in S : x^{-1}y \in N\,\, \text{and}\,\, x^0 = y^0\}$. We denote the cryptogroup $S/ \rho_{_N}$ by $S/N$ and the $x \rho_{_N}$-class by $(xN)^*$.
\end{notation}  

\begin{remark}
Let $S$ be a cryptogroup and $N$ be a full normal subcryptogroup of $S$. Then it is easy to verify that $E(S/N) = \{(xN)^* : x \in N\}$.
\end{remark}

\begin{lemma}
Let $S$ be a cryptogroup and $N$ be a full normal subcryptogroup of $S$. Then any full normal subcryptogroup of $S/N$ is of the form $M/N$, where $M$ is a full normal subcryptogroup of $S$ containing $N$.
\end{lemma}

\begin{proof}
Let $K$ be a full normal subcryptogroup of $S/N$. Let $M = \{x \in S : (xN)^* \in K\}$. Clearly, $M \neq \emptyset$. Let $x,y \in M$. Then $(xN)^*, (yN)^* \in K$. Since $K$ is a subcryptogroup of $S/N$, we have $(xyN)^*, (x^{-1}N)^* = ((xN)^*)^{-1} \in K$ and thus $xy, x^{-1} \in M$. Thus $M$ is a subcryptogroup of $S$. Now, $K$ being full, $M$ is also full. Let $s \in S, t \in M$. Then $(sN)^* \in S/N$ and $(tN)^* \in K$. Again, $K$ being normal in $S/N$, $(sN)^* (tN)^* ((sN)^*)^{-1} \in K$, i.e., $(sts^{-1}N)^* \in K$. This implies that $sts^{-1} \in M$. Therefore, $M$ is a full normal subcryptogroup of $S$. Also, for each $n \in N$, $(nN)^* \in E(S/N) \subseteq K$ implies $n \in M$. So $N \subseteq M$. Moreover, $M/N = \{(aN)^* : a \in M\} = \{(aN)^* : (aN)^* \in K\} = K$. 
\end{proof}

We now state an elementary result without proof.

\begin{lemma}
Let $S$ be a cryptogroup and $N$ be a full normal subcryptogroup of $S$. Let $H_{_x}$ and $H_{_{(xN)^*}}$ be the $\mathcal{H}$-classes containing $x$ and $(xN)^*$ in $S$ and $S/N$ respectively. Then $\psi (H_{_x}) = H_{_{(xN)^*}}$, where $\psi : \underset{x  \, \mapsto \, (xN)^*}{\overset {S \longrightarrow S/N}{}}$ is the natural epimorphism.
\end{lemma}


Let $S$ be a topological cryptogroup and $N$ be a full normal subcryptogroup of $S$. Consider the natural mapping  $\psi : \underset{x  \, \mapsto \, (xN)^*}{\overset {S \longrightarrow S/N}{}}$. Now, we define a topology on $S/N$ as : a subset $U$ of $S/N$ is open in $S/N$ if and only if $\psi^{-1}(U)$ is open in $S$. In this case, $\psi$ is a quotient mapping and $S/N$ is a quotient space.

Similar to the \cite[Theorem 5.10]{mp}, we can prove the following result.

\begin{theorem}
If $S$ is a band of topological groups and $N$ is a full normal subcryptogroup of $S$, then $S/N$ is a band of topological groups. 
\end{theorem}

\begin{theorem}
Let $S$ be a band of topological groups and $N$ be a full normal subcryptogroup of $S$. Then the following conditions are equivalent :

(i) $S/N$ is Hausdorff;

(ii) $\rho_{_N}$ is closed in $S\times S$;

(iii) $N$ is closed in $S$.
\end{theorem}

\begin{proof}
$(i) \Rightarrow (ii)$: Let $S/N$ be Hausdorff. Let $(a,b) \in \overline{\rho_{_N}}$. If possible, let $(a,b) \notin \rho_{_N}$. Then $(aN)^* \neq (bN)^*$. Since $S/N$ is Hausdorff, there exist open sets $U$ and $V$ containing $(aN)^*$ and $(bN)^*$ respectively in $S/N$ such that $U \cap V = \emptyset$. This implies $\psi^{-1}(U)$ and $\psi^{-1}(V)$ are open sets containing $a$ and $b$ respectively in $S$, where $\psi : S \longrightarrow S/N$ is the natural epimorphism. From this it follows that $\psi^{-1}(U) \times \psi^{-1}(V)$ is an open set containing $(a,b)$ in $S \times S$. Since $(a,b) \in \overline{\rho_{_N}}$, we must have $\psi^{-1}(U) \times \psi^{-1}(V) \cap \rho_{_N} \neq \emptyset$. Let $(x, y) \in \psi^{-1}(U) \times \psi^{-1}(V) \cap \rho_{_N}$. Then  $(xN)^* \in U$, $(yN)^* \in V$ and $(xN)^* = (yN)^*$. This implies $U \cap V \neq \emptyset$, a contradiction. Therefore, $(a,b) \in \rho_{_N}$ and hence $\rho_{_N}$ is closed.

$(ii) \Rightarrow (iii)$: Let $\rho_{_N}$ be closed in $S \times S$. Let $x \in \overline{N}$. If possible, let $x \notin N$. Then $(x^0)^{-1}x = x \notin N$. This implies that $(x^0, x) \notin \rho_{_N}$. Now $\rho_{_N}$ being closed, there is a basic open set $U_1 \times V_1$ containing $(x^0, x)$ in $S \times S$ such that $(U_1 \times V_1) \cap \rho_{_N} = \emptyset$. Since $\varphi : \underset{x  \, \mapsto \, x^0}{\overset {S \longrightarrow S}{}}$ is continuous, $\varphi^{-1}(U_1)$ is an open set containing $x$ in $S$. Then $\varphi^{-1}(U_1) \cap V_1$ is an open set containing $x$ in $S$. Now $x \in \overline{N}$ implies that $(\varphi^{-1}(U_{_1}) \cap V_{_1}) \cap N \neq \emptyset$. Let $p \in (\varphi^{-1}(U_{_1}) \cap V_{_1}) \cap N$. Then $p^0 \in U_{_1}$, $p \in V_1 \cap N$ and hence  $(p^0, p) \in (U_1 \times V_1) \cap \rho_{_N}$, a contradiction. Hence $x \in N$ and so $N$ is closed. 

$(iii) \Rightarrow (i)$: Let $N$ be closed. Let $(uN)^*, (vN)^* \in S/N$ with $(uN)^* \neq (vN)^*$. Then $(u,v) \notin \rho_{_N}$. This implies either $u^{-1}v \notin N$ or $u^0 \neq v^0$.  If $u^0 \neq v^0$, then $H_{_u}$ and $H_{_v}$ are disjoint open sets containing $u$ and $v$ respectively in $S$. Now, the natural mapping $\psi : S \longrightarrow S/N$ being open, $\psi (H_{_u})$ and $\psi (H_{_v})$ are disjoint open  sets containing $(uN)^*$ and $(vN)^*$ respectively in $S/N$. Hence $S/N$ is Hausdorff. On the other hand, if $u^{-1}v \notin N$, then $N$ being closed, there is an open set $G$ containing $u^{-1}v$ in $S$ such that $G \cap N = \emptyset$. Now, since $S$ is a topological cryptogroup, there exist open sets $U_2$ and $V_2$ containing $u$ and $v$ respectively in $S$ such that $U_2^{-1}V_2 \subseteq G$. Then $\psi (U_2)$ and $\psi (V_2)$ are open sets containing $(uN)^*$ and $(vN)^*$ respectively in $S/N$. If possible, let $\psi (U_2) \cap \psi (V_2) \neq \emptyset$ and let $(xN)^* \in \psi (U_2) \cap \psi (V_2)$. Then $(xN)^* = (pN)^*$ and $(xN)^* = (qN)^*$, for some $p \in U_2$ and $q \in V_2$. From this, we have $p^{-1}q \in U_2^{-1}V_2 \subseteq G$ and $p^{-1}q \in N$. This implies $G \cap N \neq \emptyset$, a contradiction. Hence $\psi (U_2)$ and $\psi (V_2)$ are disjoint open sets containing $(uN)^*$ and $(vN)^*$ respectively in 
$S/N$. Consequently, $S/N$ is Hausdorff. 
\end{proof}

\noindent \textbf{Acknowledgement:}

\noindent The research of the second author was supported by the Council for Scientific and Industrial Research, India. CSIR Award no. : 09/028(1001)/2017-EMR-1.

\end{document}